\documentclass{article}

\usepackage{amsmath,amssymb,epsf,wrapfig}
\usepackage[dvipdfm]{graphicx}
\def\qed{\hfill$\square$}
\newtheorem{theorem}{Theorem}
\newtheorem{corollary}[theorem]{Corollary}
\newtheorem{lemma}[theorem]{Lemma}
\newtheorem{proposition}[theorem]{Proposition}

\newtheorem{remark}[theorem]{Remark}

\newtheorem{definition}[theorem]{Definition}

\begin{document}

\title{Chart description for genus-two Lefschetz fibrations and a theorem on their stabilization}
\author{Seiichi Kamada\\
Department of Mathematics, Hiroshima University, \\
Higashi-Hiroshima, Hiroshima 739-8526, Japan}

\date{}

\maketitle

\begin{abstract}
Chart descriptions are a 
graphic method to describe monodromy representations of various topological objects.  Here we introduce a chart description for genus-two Lefschetz fibrations, and show that any genus-two Lefschetz fibration can be stabilized by fiber-sum with certain basic Lefschetz fibrations. 
\end{abstract}

\section{Introduction}

Chart descriptions were originally introduced in order to describe $2$-dimensional braids  in \cite{Kam92, Kam96} (cf. \cite{Kam02}).  
In \cite{KMMW}, a chart description for genus-one Lefschetz fibrations was introduced and an elementary proof of Matsumoto's  classification theorem was given.  
At the third JAMEX meeting in Oaxaca, Mexico, 2004, the author  generalized it to a method 
describing any monodromy representation \cite{Kam07}.  
Here we introduce a chart description for genus-two Lefschetz fibrations, and show that any genus-two Lefschetz fibration can be stabilized by fiber-sum with certain basic Lefschetz fibrations.  Our result was partially announced at `The Second East Asia School of Knots and 
Related Topics in Geometric Topology' in Dalian, China, 2005.  

\section{Lefschetz fibrations} 

Let $M$ and $B$ be compact, connected, and 
oriented smooth $4$-manifold and $2$-manifold, respectively.   
Let $f : M \to B$ be a smooth map with $\partial M = f^{-1}(\partial B)$.    
A critical  point $p$ is called a {\it Lefschetz singular point} of 
{\it positive type} (or of {\it negative type}, respectively) 
if there exist local complex 
coordinates $z_1, z_2$ around $p$ and 
a local complex coordinate $\xi$ around $f(p)$  
such that $f$ is locally written as $\xi= f(z_1,z_2) =z_1z_2$ (or $\overline{z_1}z_2$, resp.).  
We call  $f$ 
a  (smooth or differentiable) {\it  Lefschetz fibration\/} if 
all critical points are Lefschetz singular points and if 
there exists exactly one critical point in the preimage of each critical value.  

A {\it general fiber\/} is the preimage of a regular value of $f$.  
A {\it singular fiber\/} of {\it positive type} (or  {\it negative type}, resp.) is the preimage 
of a critical value which contains a Lefschetz singular point of positive type (or negative type, resp.).  
A singular fiber is obtained by shrinking a simple loop, called a vanishing cycle, on a general fiber.  
In this paper we assume that a Lefschetz fibration is `{relatively minimal}', i.e., all vanishing cycles are essential loops.  
We say that a singular fiber is of  {\it type {\rm I}} or of {\it type {\rm II}} \/ if the vanishing cycle is a 
non-separating loop or a separating loop, respectively.  

A singular fiber is of {\it type ${\rm I}^+$} if it is of type I and of positive type.  Similarly 
{\it type ${\rm I}^-$},  {\it type ${\rm II}^+$}  {\it type ${\rm II}^-$} are defined.   
We denote by 
$n_{\rm I}^{+}(f)$, $n_{\rm I}^{-}(f)$, $n_{\rm II}^{+}(f)$, and $n_{\rm II}^{-}(f)$, 
the numbers of singular fibers of $f$ of type ${\rm I}^+$, ${\rm I}^-$, ${\rm II}^+$, and ${\rm II}^-$, respectively.  
A Lefschetz fibration is called  {\it irreducible} if 
every singular fiber is of type I, i.e.,  $n_{\rm II}^+(f) = n_{\rm II}^-(f) = 0$.   
A Lefschetz fibration is called {\it chiral} or {\it symplectic} if every singular fiber is of positive type, i.e.,  $n_{\rm I}^-(f) = n_{\rm II}^-(f)=0.$

Let $f : M \to B$ be a Lefschetz fibration, and  
$\Delta =\{q_1, \dots, q_n\}$ the set of critical values.  
Let $\rho: \pi_1(B \setminus \Delta, q_0) \to MC$ be the monodromy representation of $f$, where $q_0$ is a base point of $B \setminus \Delta$ and $MC$ is the mapping class group of the fiber $f^{-1}(q_0)$.  
Consider a Hurwitz arc system for 
$\Delta$, say ${\cal A}= (A_1, \dots, A_n)$; each $A_i$ is an embedded arc in $B$ connecting $q_0$ 
and a point of $\Delta$ such that $A_i \cap A_j= \{q_0\}$ for $i \ne j$, and they appear in this order around $q_0$.  
When $B$ is a $2$-sphere or a $2$-disk, the system ${\cal A}$ determines a system of generators of $\pi_1(B \setminus \Delta, q_0)$, say $(a_1, \dots, a_n)$.  
We call $( \rho(a_1), \dots, \rho(a_n) )$ a {\it Hurwitz system} of $f$.  
For details on Hurwitz systems, refer to  \cite{Auroux2003, GS, Matsu96, Moi81, ST03}, etc.

\section{Main result} 

\begin{figure}[h]
\begin{center}
\mbox{\epsfxsize=5cm \epsfbox{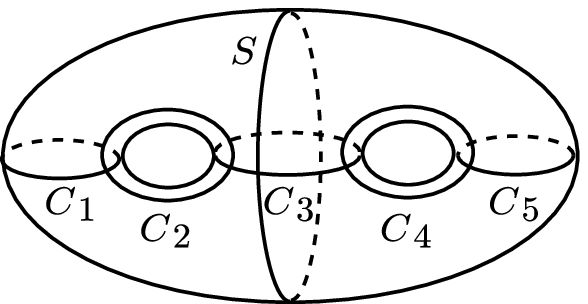}}
\end{center}
\vspace{-0.5cm}
\caption{}
\label{sfg01}
\end{figure}

Let $\zeta_i$ $(i=1,\ldots, 5)$ be positive Dehn twists along 
the loops $C_i$ $(i=1,\ldots, 5)$ illustrated in Figure~\ref{sfg01}.  
The mapping class group $MC$ of a genus-two Riemann surface is generated by 
$\zeta_1, \zeta_2, \zeta_3, \zeta_4, \zeta_5$, and 
the following relations are defining relations (cf. \cite{Birman}). 
\begin{eqnarray}
&& \zeta_i \zeta_j = \zeta_j \zeta_i     \quad \mbox{    if $|i-j|\geq 2$, } \label{eq:01} \\
&& \zeta_i \zeta_{i+1} \zeta_i = \zeta_{i+1} \zeta_i \zeta_{i+1}    \quad \mbox{    for $i=1, \dots, 4$, } \label{eq:02} \\
&& \iota^2 =1    \quad \mbox{  where $\iota= \zeta_1 \zeta_2 \zeta_3 \zeta_4 \zeta_5^2 \zeta_4 \zeta_3 \zeta_2 \zeta_1$,}\label{eq:03} \\
&&(\zeta_1 \zeta_2 \zeta_3 \zeta_4 \zeta_5)^6=1, \label{eq:04} \\
&&\iota  \,  \zeta_i =
 \zeta_i \,  \iota  
\quad \mbox{    for $i=1, \dots, 5$.} \label{eq:05} 
\end{eqnarray}

Let $\sigma$ be a positive Dehn twist along the loop $S$ illustrated in Figure~\ref{sfg01}.  
Then $
\sigma= (\zeta_1 \zeta_2)^6$. 

If $(g_1, \ldots, g_n)$ is a Hurwitz system of a genus-two Lefschetz fibration, then 
each $g_j$ is a conjugate of $\zeta_i$ or $\zeta_i^{-1}$, or a conjugate of $\sigma$ or $\sigma^{-1}$.  

Now we define basic Lefschetz fibrations.  

\begin{definition}[cf. \cite{Auroux2003, Auroux2005, Matsu96, ST03}]{\rm 
{\it Basic Lefschetz fibrations},  $f_0$, $f_1$, $f_2$, $f'_1$ and $f'_2$,  are  
Lefschetz fibrations over $S^2$ whose Hurwitz systems are 
\begin{itemize}
\item[(1)] 
 $ W_0 = (T)^2 $ where $T = (\zeta_1, \zeta_2, \zeta_3, \zeta_4, \zeta_5, 
\zeta_5, \zeta_4, \zeta_3, \zeta_2, \zeta_1)$, 

\item[(2)]  $ W_1 = (\zeta_1, \zeta_2, \zeta_3, \zeta_4, \zeta_5)^6$, 

\item[(3)] $ W_2 = (\sigma, (\zeta_3, \zeta_4, \zeta_5, \zeta_2, \zeta_3, \zeta_4, \zeta_1, \zeta_2, \zeta_3)^2, T) $,  

\item[(4)]  $W'_1 = (\zeta_1, \zeta_1^{-1})$,  

\item[(5)]  $W'_2 = (\sigma, \sigma^{-1})$, 
\end{itemize}
respectively.  
}\end{definition}

For example, $f_0$ has $20$ singular fibers, which are of type ${\rm I}^+$.  Thus $f_0$ is chiral and irreducible.   

\newcommand{\lw}[1]{\smash{\lower2.0ex\hbox{#1}}}
\begin{center}
\renewcommand{\arraystretch}{1.2}
 \begin{tabular}{|l@{\quad\vrule width0.8pt}c|c|c|c|c|c|}
 \hline 
 \lw{LF} & \multicolumn{4}{c|}{\# of sing. fib.} & \lw{chiral} & \lw{irreducible}   \\
 \cline{2-5}  & \, $n_{\rm I}^+$ & $n_{\rm I}^-$ & $n_{\rm II}^+$ & $n_{\rm II}^-$ & &\\
\noalign{\hrule height 0.8pt} 
 $f_0$ & 20 & 0 & 0 & 0 & $\bigcirc$ & $\bigcirc$ \\
 $f_1$ & 30 & 0 & 0 & 0 & $\bigcirc$ & $\bigcirc$ \\
 $f_2$ & 28 & 0 & 1 & 0 & $\bigcirc$ & $\times$ \\
 $f'_1$ & 1 & 1 & 0 &  0 & $\times$ & $\bigcirc$ \\
 $f'_2$ & 0 & 0 & 1 & 1 & $\times$ & $\times$ \\
\hline 
\end{tabular}
\end{center}

For two Lefschetz fibrations $f$ and $f'$ over $S^2$, we denote by $f \# f'$ the fiber-sum of $f$ and $f'$.  
By ${\#} m f$ for a positive integer $m$, we mean the fiber-sum of $m$ copies of $f$. 

\begin{theorem}\label{thm:main}
Let $f$ be a genus-two Lefschetz fibration over $S^2$.  
Suppose that $n_{\rm II}^+(f) \geq n_{\rm II}^-(f)$. Then 
\begin{itemize}
\item[{\rm (1)}] 
${\cal E}(f) := n_{\rm I}^+(f) - n_{\rm I}^-(f) -28 (n_{\rm II}^+(f) - n_{\rm II}^-(f))$  
is a multiple of $10$.  
\end{itemize}  
(Hence we can define  the  {\it parity}, $\epsilon(f) \in \{0,1\}$, by 
$\epsilon(f) \equiv  {\cal E}(f)/10 \mod{2}$. )   

\begin{itemize}
\item[{\rm (2)}] 
There exists a positive integer $m_0$ such that  
for any integer $m \geq m_0$, 
$$ 
f \,\# \, m   \,  f_0 
 \cong 
       {\#}  \, (a+m)  \,  f_0 
 \,  \# \,  b  \,  f_1 
 \,  \#  \, c   \,  f_2
 \,  \#  \, d   \, f'_1
 \,  \#  \,  e  \,  f'_2 $$
for some non-negative integers $a, b, c, d$ and  $e$. 
\item[{\rm (3)}] 
In $(2)$,  it holds that 
 $c = n_{\rm II}^+(f) - n_{\rm II}^-(f)$,  
$d = n_{\rm I}^-(f)$ and $e = n_{\rm II}^-(f)$.  
Although $a$ and $b$ are not determined uniquely,  we can take 
$b = \epsilon(f) \in \{0,1\}$, 
and then 
$a = ({\cal E}(f) -30 \epsilon(f))/20$. 
\item[{\rm (4)}] 
If $n_{\rm II}^-(f)=0$, then we may take $m_0$ in $(2)$ to be $n_{\rm I}^-(f) + 2   n_{\rm II}^+(f)  +1$.  
\end{itemize}
\end{theorem}

 \begin{remark}{\rm 
If $f$ is chiral and irreducible, then 
$n_{\rm I}^-(f) = n_{\rm II}^+(f) = n_{\rm II}^-(f) = 0$  and  by (4) 
we may assume $m_0 = 1$. 
Thus,  we have 
$$f  \, \#    \,  f_0 
\cong 
        \, {\#}   \,  (a+1)   \,  f_0
  \,  \#   \,  b   \,   f_1. $$
This is due to B. Siebert and G. Tian \cite{ST03}.  Our proof concerning the assertion (4) of Theorem~\ref{thm:main} is based on their result.  In Section~\ref{sect:concluding},  we observe that Theorem~\ref{thm:main} except the assertion (4) can be proved without Siebert and Tian's result.  

If $f$ is chiral, then 
$n_{\rm I}^-(f) = n_{\rm II}^-(f) = 0$. 
By Theorem~\ref{thm:main}, we have \\ 
$$f    \,  \#   \,  m   \,  f_0  
\cong 
        \,  {\#}   \,  (a+m)   \,  f_0
  \,  {\#}   \,  b   \,  f_1 
  \,  {\#}   \,  c   \, f_2. $$ 
This is due to D. Auroux \cite{Auroux2003}.   Here  $m$ is any integer with $m \geq 2 n_{\rm II}^+(f)  +1$.  
}\end{remark}


\section{Chart description} 

In this section we introduce a chart description for genus-two Lefschetz fibrations.  
We use the terminologies on chart description in \cite{Kam07}.  
For simplicity's sake, we only consider 
genus-two Lefschetz fibrations over $B$ such that $\partial B$ is empty or connected, and if  $\partial B$ is not empty, 
we assume that the monodromy along $\partial B$ is trivial.  Unless otherwise stated, genus-two Lefschetz fibrations over $B$ are assumed to be so.  

\begin{figure}[h]
\begin{center}
\mbox{\epsfxsize=10cm \epsfbox{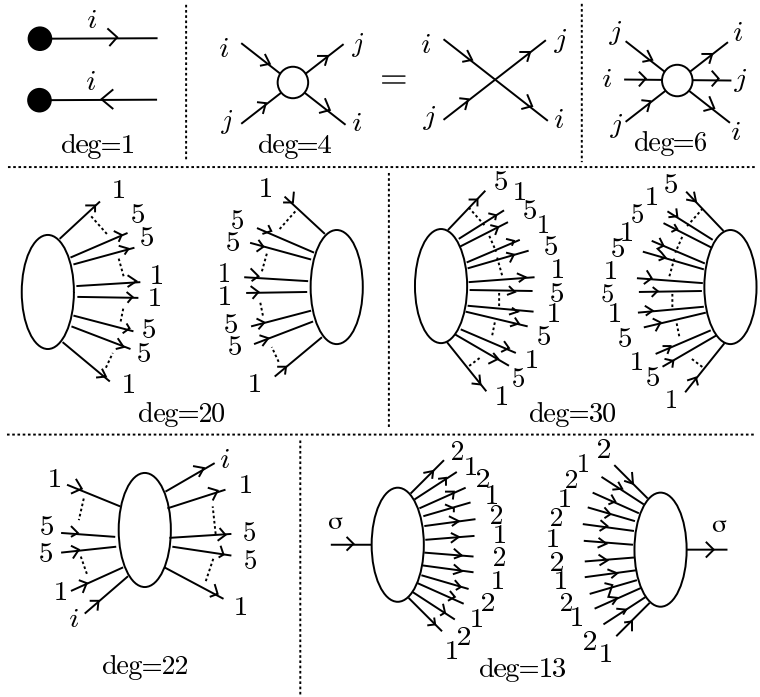}}
\end{center}
\vspace{-0.5cm}
\caption{Vertices of a chart}
\label{sfg02}
\end{figure}

\begin{definition}[cf. \cite{Kam02, Kam07, KMMW}]\label{def:chart} {\rm 
A {\it chart} in $B$ is a finite graph $\Gamma$ in $B$ (possibly being empty or having {\it hoops} that are closed edges without vertices) whose edges are labeled with an element of $\{1, 2, 3, 4, 5, \sigma\}$, and oriented so that the following conditions are satisfied (see Figure~\ref{sfg02}):   
\begin{itemize}

\item[(1)] The degree of each vertex is $1, 4, 6, 20, 30, 22$ or $13$.  
\item[(2)] For a degree-$1$ vertex, the adjacent edge is oriented outward or inward.  
\item[(3)] For a degree-$4$ vertex, two edges in each diagonal position have the same label and are oriented coherently; and the labels $i$ and $j$ of the diagonals are in $\{1, \dots, 5\}$ with $| i - j | >1$.  
\item[(4)] For a degree-$6$ vertex, the six edges are alternately labeled $i$ and $j$  in $\{1, \dots, 5\}$ with  $|i-j|=1$; and three consecutive edges are oriented outward while the other three are oriented inward.  
\item[(5)] For a degree-$20$ vertex,  the edges are labeled with  $(1,2,3,4,5,5,4,3,2,1)^2$; and all edges are oriented outward or all edges are oriented inward. 
\item[(6)] For a degree-$30$ vertex, the edges are labeled with $(1,2,3,4,5)^6$ in a counterclockwise direction 
(or clockwise direction, resp.); and all edges are oriented outward (or inward, resp.). 
\item[(7)] For a degree-$22$ vertex, the edges are labeled with $(1,2,3,4,5,5,4,3,2,1,i)^2$ in a counterclockwise direction where $i \in \{1, \dots, 5\}$; and the first $11$ edges are oriented outward and the latter ones are oriented inward.   
\item[(8)] For a degree-$13$ vertex, the edges are labeled with $((1,2)^6, \sigma)$  
in a counterclockwise direction 
(or clockwise direction, resp.); and the edges with labels $1$ and $2$ are oriented outward (or inward, resp.), and 
the edge with label $\sigma$ is oriented inward (or outward, resp.).  
\item[(9)] $\Gamma \cap \partial B = \emptyset$.   
\item[(10)] $\Gamma$ misses the base point $q_0 \in B$.  
\end{itemize}
}\end{definition}

\begin{remark}{\rm 
When we would treat genus-two Lefschetz fibrations over $B$ with $\partial B \neq \emptyset$ such that the monodromies along $\partial B$ are not  trivial, the condition (9) should be removed.  See \cite{Kam07}.  
}\end{remark}

We call a degree-$1$ vertex a {\it black vertex}.  
We say that a chart  is  {\it chiral} if every black vertex has an adjacent edge oriented outward.  
We say that a chart  is {\it irreducible} if there exist no edges with label $\sigma$.

For a chart $\Gamma$, let $\Delta_\Gamma$ be the set of black vertices.  
A chart $\Gamma$ determines a homomorphism $\pi_1(B \setminus \Delta_\Gamma, q_0) \to MC$ as in \cite{Kam07}.  
By Theorem~5 of \cite{Kam07}, we have the following theorem. 

\begin{theorem}\label{thm:chartdescription}
Let $f $ be a genus-two Lefschetz fibration  over $B$, and let $\rho$ be the monodromy representation.    Then there is a chart $\Gamma$ in $B$ such that 
the monodromy representation $\rho$  equals 
the homomorphism  $\rho_\Gamma$ determined by $\Gamma$.    
\end{theorem}

A chart $\Gamma$ as in Theorem~\ref{thm:chartdescription} is called a {\it chart description} of $f$ or a {\it chart describing } $f$.  A chart $\Gamma$ in $D^2$ is also regarded as a chart in $S^2$ in the trivial way. 

\vspace{0.5cm}

We introduce some local moves on chart descriptions.  

(C1) For a chart $\Gamma$, suppose that there exists a chart $\Gamma'$ and an embedded 2-disk, say $E$,  in $B$ such that (i) $\partial E$ intersects with $\Gamma$ and $\Gamma'$ transversely  (or do not intersect with them) avoiding their vertices, (ii) $\Gamma$ and $\Gamma'$ have no black vertices in $E$, and (iii) $\Gamma$ and $\Gamma'$ are identical outside of $E$.  Then we say that $\Gamma'$ is obtained from $\Gamma$ by a C1-move.  

\begin{figure}[h]
\begin{center}
\mbox{\epsfxsize=12cm \epsfbox{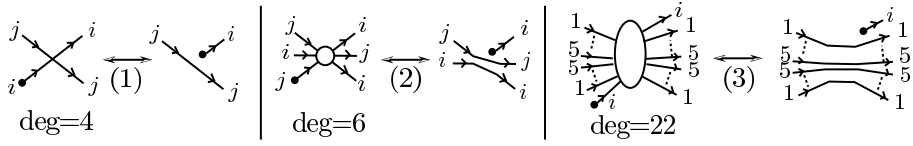}}
\end{center}
\vspace{-0.5cm}
\caption{Some chart moves}
\label{sfg03}
\end{figure}

(C2) For a chart, suppose that there is an edge $e$ joining a  degree-$4$ vertex and  a black vertex.  Remove the edge $e$ as in Figure~\ref{sfg03}(1).  
We call this local move  a C2-move.  

(C3) For a chart, suppose that  there is an edge $e$ joining a degree-$6$ vertex and a black vertex.  Suppose that $e$ is neither the middle of three edges oriented outward nor the middle of the three edges oriented inward.  Then, remove the edge as in 
Figure~\ref{sfg03}(2).  
We call this local move  a C3-move.  

(C4) In a chart, suppose that  there is an edge $e$ joining a degree-$22$ vertex and a black vertex.  Suppose that $e$ is one of the two edges labeled $i$ in Figure~\ref{sfg02}.  Then, remove the edge as in 
Figure~\ref{sfg03}(3).  
We call this local move  a C4-move.  

\vspace{0.5cm}

When $\partial B \neq \emptyset$ and the base point $q_0$ is in $\partial B$, we introduce another move. 

(C5) Suppose that $\partial B \neq \emptyset$ and $q_0 \in \partial B$.  Let $\Gamma'$ be a chart that is the union of  a chart $\Gamma$ and some hoops which are parallel to and sufficiently near $\partial B$.  Then we say that $\Gamma'$ is obtained from $\Gamma$ by a C5-move.

\begin{definition}\label{def:chartmove}{\rm 
(1) {\it Chart moves} are C1-moves, C2-moves, C3-moves, C4-moves and their inverse moves.  

(2) Two charts in $B$ are said to be {\it chart move equivalent} (with respect to the base point $q_0$) if they are related by a finite sequence of chart moves and ambient isotopies of $B$ rel $q_0$, where we assume that chart moves are applied in embedded 2-disks in $B$ missing $q_0$. 

(3)   Two charts in $B$ are said to be {\it chart move equivalent up to conjugation} (with respect to the base point $q_0$) if they are related by a finite sequence of chart moves, C5-moves and ambient isotopies of $B$ rel $q_0$.  (It is not necessary to assume that chart moves are applied in embedded 2-disks in $B$ missing $q_0$.)
}\end{definition}

We say that two monodromy representations $\rho: \pi_1(B \setminus \Delta, q_0) \to MC$ and $\rho': \pi_1(B \setminus \Delta', q_0) \to MC$ are {\it equivalent} if there is a diffeomorphism $h : (B, q_0) \to (B, q_0)$ which is isotopic to the identity map rel $q_0$ such that $h(\Delta)= \Delta'$ and 
$\rho = \rho' \circ h_\ast$, where $h_\ast : \pi_1(B \setminus \Delta, q_0) \to  \pi_1(B \setminus \Delta', q_0)$ is the induced isomorphism.  

We say that two monodromy representations $\rho: \pi_1(B \setminus \Delta, q_0) \to MC$ and $\rho': \pi_1(B \setminus \Delta', q_0) \to MC$ are {\it equivalent up to conjugation} if there is an inner-automorphism of $MC$, say $t$, and 
there is a diffeomorphism $h : (B, q_0) \to (B, q_0)$ which is isotopic to the identity map rel $q_0$ such that $h(\Delta)= \Delta'$ and 
$\rho = t \circ \rho' \circ h_\ast$.  

\vspace{0.5cm}

C1-moves in this paper are called  {\it chart moves of type $W$} in Definition~7 of \cite{Kam07}.  
C2-moves, C3-moves, C4-moves, C5-moves are not given explicitly in \cite{Kam07}.  
However, as shown in Fig. 22 and 23 of \cite{Kam07}, C2-moves and  C3-moves are 
equivalent to some 
local moves called {\it chart moves of transition} in Definition~14 of \cite{Kam07}.  
C4-moves are also equivalent to chart moves of transition in the sense of \cite{Kam07}.  
Thus, as stated in Section~8 of \cite{Kam07},  we see that if two charts are chart move equivalent in our sense (Definition~\ref{def:chartmove} (2)) then the monodromy representations determined by them are equivalent.  
C5-moves are equivalent to {\it chart moves of conjugacy} in (3) and (4) of Fig.~17 of \cite{Kam07}.  Again as in Section~7 of \cite{Kam07},  we see that if two charts are chart move equivalent up to conjugation (Definition~\ref{def:chartmove} (3)) then the monodromy representations determined by them are equivalent up to conjugation.  

Thus we have the following.  

\begin{theorem}
For two charts in $B$, if they are chart move equivalent (or chart move equivalent up to conjugation, resp.) then 
the monodromy representations determined by them are equivalent (or equivalent up to conjugation, resp.),  
and hence the Lefschetz fibrations described by them are isomorphic.  
\end{theorem}

\begin{remark}{\rm 
By Theorem~16 of \cite{Kam07}, we see that two charts determine equivalent monodromy representations  if and only if they are related by C1-moves (chart move of type $W$), chart moves of transition, and ambient isotopies of $B$ rel $q_0$.  It is unknown to  the author whether all chart moves of transition are consequence of our chart moves.  
}\end{remark}

We say that a black vertex of a chart $\Gamma$ is of {\it type ${\rm I}^+$}, 
{\it type ${\rm I}^-$},  {\it type ${\rm II}^+$} or  {\it type ${\rm II}^-$} if 
the adjacent edge is labeled in $\{1, \dots, 5\}$ and oriented outward, if 
the adjacent edge is labeled in $\{1, \dots, 5\}$ and oriented inward, if 
the adjacent edge is labeled $\sigma$ and oriented outward, or if 
the adjacent edge is labeled $\sigma$ and oriented inward, respectively.   

When $\Gamma$ is a chart description of a genus-two Lefschetz fibration $f : M \to B$,  
black vertices correspond to critical values of $f$, and 
the types of the vertices are the same with the types of the singular fibers over the corresponding critical values.   
For a chart $\Gamma$,  
we denote by 
$n_{\rm I}^{+}(\Gamma)$, $n_{\rm I}^{-}(\Gamma)$, $n_{\rm II}^{+}(\Gamma)$, and $n_{\rm II}^{-}(\Gamma)$, 
the numbers of black vertices of type ${\rm I}^+$, 
type ${\rm I}^-$,  type ${\rm II}^+$ and type ${\rm II}^-$, respectively.  They are equal to 
$n_{\rm I}^{+}(f)$, $n_{\rm I}^{-}(f)$, $n_{\rm II}^{+}(f)$, and $n_{\rm II}^{-}(f)$, respectively. 

\vspace{0.5cm}
If a chart $\Gamma$ is irreducible, then it is obvious that $n_{\rm II}^{+}(\Gamma) = n_{\rm II}^{-}(\Gamma) =0$.  The converse is not true.  However we have the following.  

\begin{lemma}\label{lem:irreduciblechart}
Every chart $\Gamma$ with $n_{\rm II}^{+}(\Gamma) = n_{\rm II}^{-}(\Gamma) =0$ is chart move equivalent to an irreducible chart.  
\end{lemma}

\begin{figure}[h]
\begin{center}
\mbox{\epsfxsize=12cm \epsfbox{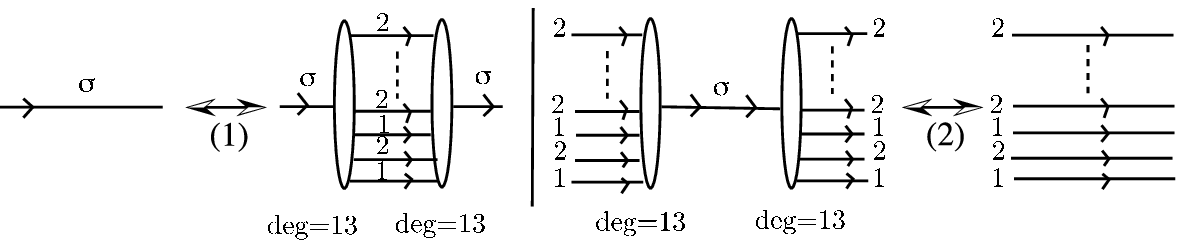}}
\end{center}
\vspace{-0.5cm}
\caption{}
\label{sfg04}
\end{figure}

{\it Proof.}  We can replace every hoop labeled $\sigma$ into  $12$ parallel hoops with labels $1$ or $2$ by 
a chart move depicted in Figure~\ref{sfg04} (1) followed by one  in Figure~\ref{sfg04} (2).  
Every edge labeled $\sigma$ whose endpoints are degree-$13$ vertices is also removed by the latter move.  \qed  

\begin{proposition}\label{prop:chartdescriptionchiralirreducible}
A chiral (or irreducible, resp.) genus-two Lefschetz fibration has a chart description which is chiral 
(or irreducible, resp.).  
\end{proposition}

{\it Proof.} If $f$ is chiral, local monodromies around the critical values are all positive Dehn twists.  By the definition of a chart description, the adjacent edges of the black vertices are oriented outward.  Thus any chart description of $f$ is chiral.  If $f$ is 
irreducible, local monodromies around the critical values are Dehn twists along non-separating simple loops, which are conjugates of $\zeta_1, \dots, \zeta_5$ and their inverses.  Thus any chart description $\Gamma$ of $f$ satisfies $n_{\rm II}^{+}(\Gamma) = n_{\rm II}^{-}(\Gamma) =0$.  By Lemma~\ref{lem:irreduciblechart}, it changes to an irreducible one. \qed

\begin{figure}[h]
\begin{center}
\mbox{\epsfxsize=9cm \epsfbox{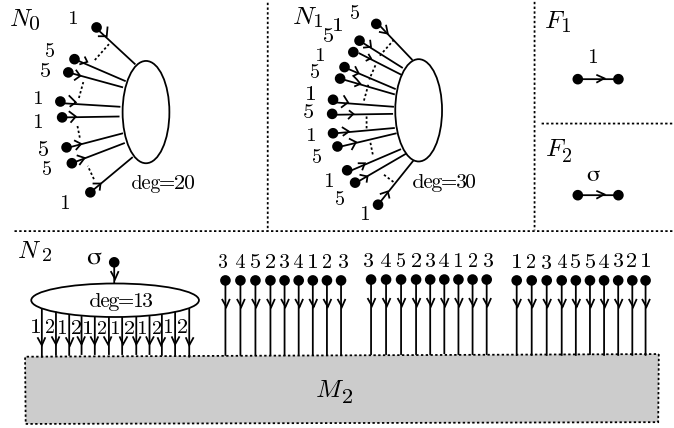}}
\end{center}
\vspace{-0.5cm}
\caption{Charts $N_0$, $N_1$, $N_2$, $F_1$ and $F_2$ describing 
$f_0$, $f_1$, $f_2$, $f'_1$ and $f'_2$}
\label{sfg05}
\end{figure}

In Figure~\ref{sfg05}, we show charts $N_0$, $N_1$, $N_2$, $F_1$ and $F_2$ describing 
$f_0$, $f_1$, $f_2$, $f'_1$ and $f'_2$.  
We call $N_0$ a (positive) {\it nucleon of degree-$20$} and $N_1$ a (positive) {\it nucleon of degree-$30$}.  
The region named $M_2$ is an arbitrary chart consisting of edges with labels in $\{1, \dots, 5\}$ and vertices whose degrees are in $\{4, 6, 20, 22\}$.  
(There exists such a chart $M_2$,  Lemma~\ref{lem:change}.)  
A {\it free edge} means a chart consisting two black vertices and a single edge connecting them.  
$F_1$ and $F_2$ are free edges.  

\vspace{0.5cm}

Let $\Gamma$ and $\Gamma'$ be charts in $B=D^2$.  Divide $D^2$ into $2$-disks $D^2_1$ and $D^2_2$ by a properly embedded arc in $D^2$.  Put a small copy of $\Gamma$ in $D^2_1$ and a small copy of $\Gamma'$ in $D^2_2$.  We have a new chart in $D^2 = D^2_1 \cup D^2_2$.  We call it the {\it product} of $\Gamma$ and $\Gamma'$ and denote it by $\Gamma \oplus \Gamma'$.  
We say that $\Gamma$ is a {\it factor} of $\Gamma \oplus \Gamma'$.  
The chart $\Gamma \oplus \Gamma'$ is a chart description of the fiber sum $f \# f'$ of the Lefschetz fibreations $f$ and $f'$ described by $\Gamma$ and $\Gamma'$.   We denote by $n \Gamma$ the product $\Gamma \oplus \cdots \oplus \Gamma$ of $n$ copies of $\Gamma$.  (When $B=D^2$, the fiber sum $f \# f'$ of $f$ and $f'$ over $B$ is defined by using the boundary connected sum of the base spaces.)

\begin{theorem}\label{thm:chart1}
Let $\Gamma$ be a chart in $B=D^2$.  Suppose that $n_{\rm II}^+(\Gamma) \geq n_{\rm II}^-(\Gamma)$. 
Then there exists a positive integer $m_0$ such that for any integer $m \geq m_0$, 
the chart $\Gamma \oplus m  \, N_0$ is chart move equivalent to 
$$\Gamma'  \oplus (n_{\rm II}^+(\Gamma)-  n_{\rm II}^-(\Gamma)) N_2  \oplus n_{\rm I}^-(\Gamma) \, F_1  
\oplus  n_{\rm II}^-(\Gamma) F_2  $$
for some chart $\Gamma'$ with $n_{\rm I}^-(\Gamma')=n_{\rm II}^+(\Gamma') = n_{\rm II}^-(\Gamma')=0$ 
such that $\Gamma'$ has $N_0$ as a factor.  
Moreover if $n_{\rm II}^-(\Gamma)=0$, we may take $m_0$ to be $n_{\rm I}^-(\Gamma)  +   2 n_{\rm II}^+(\Gamma) +1 $.  
\end{theorem}

We prove Theorem~\ref{thm:chart1} in Section~\ref{sect:chart1proof}.  

\begin{corollary}\label{cor:chart1}
Let $f$ be a genus-two Lefschetz fibration over $B=D^2$ (or $S^2$) with $n_{\rm II}^+(f) \geq n_{\rm II}^-(f)$. 
Then there exists a positive integer $m_0$ such that for any integer $m \geq m_0$, 
the fiber sum $f  \# m  f_0$ is equivalent to 
$$f'  \# (n_{\rm II}^+(f)-  n_{\rm II}^-(f)) f_2  \# n_{\rm I}^-(f) \, f'_1  
\#  n_{\rm II}^-(f) f'_2  $$ 
for some chiral and irreducible genus-two Lefschetz fibration $f'$ over $B=D^2$ (or $S^2$) 
such that the monodromy representation of $f'$ is transitive.  
Moreover if $n_{\rm II}^-(f)=0$, we may take $m_0$ to be $n_{\rm I}^-(f)  +   2 n_{\rm II}^+(f) +1 $. 
\end{corollary}

\begin{lemma}\label{lem:change}
There is a chart satisfying the condition of $M_2$.  
\end{lemma}

\begin{figure}[h]
\begin{center}
\mbox{\epsfxsize=9.0cm \epsfbox{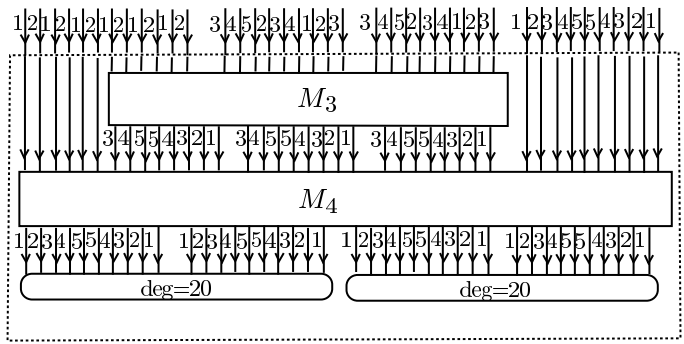}}
\end{center}
\vspace{-0.5cm}
\caption{Chart $M_2$}
\label{sfg06}
\end{figure}

\begin{figure}[h]
\begin{center}
\mbox{\epsfxsize=6.5cm \epsfbox{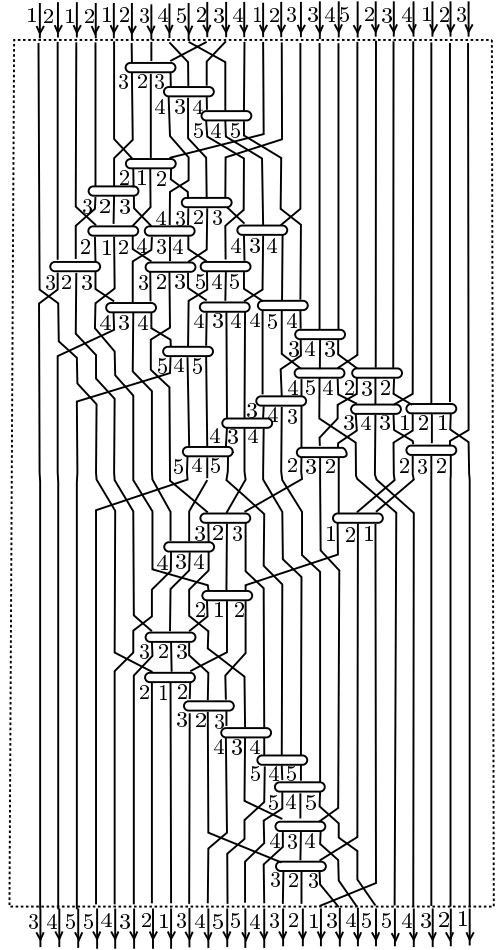}}
\end{center}
\vspace{-0.5cm}
\caption{Chart $M_3$}
\label{sfg07}
\end{figure}

\begin{figure}[h]
\begin{center}
\mbox{\epsfxsize=9.0cm \epsfbox{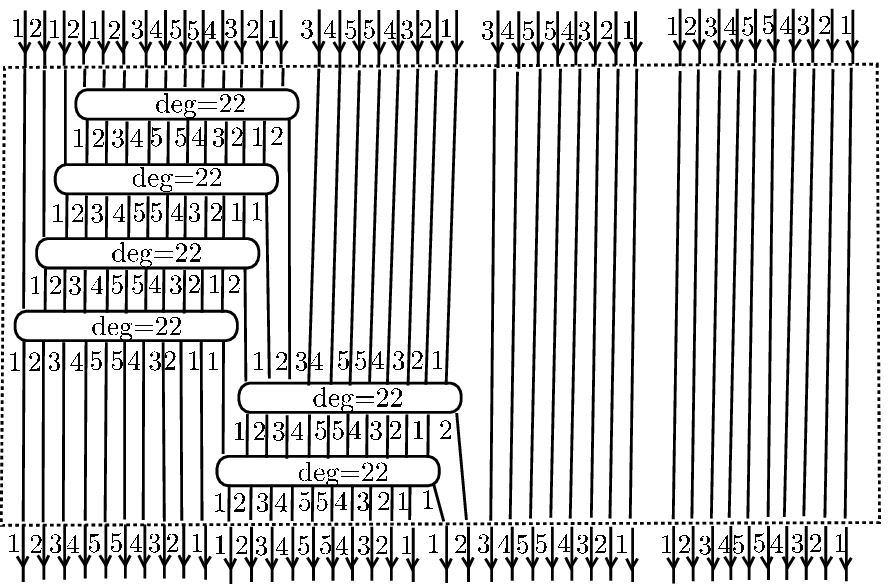}}
\end{center}
\vspace{-0.5cm}
\caption{Chart $M_4$}
\label{sfg08}
\end{figure}

{\it Proof.} 
See Figure~\ref{sfg06} where $M_3$ and $M_4$ are charts depicted  in 
Figures~\ref{sfg07} and~ \ref{sfg08}.  
\qed 

\section{Proof of Theorem~\ref{thm:chart1}}\label{sect:chart1proof}

\begin{definition}{\rm 
A chart $\Gamma$ in a $2$-disk   is {\it nomadic  with respect to} a chart $\Gamma_0$ in $B$ if 
for any two regions of the complement $B \setminus \Gamma_0$, say $R_1$ and $R_2$,  
the chart $\Gamma_0$ together with a small copy of $\Gamma$ in $R_1$ is chart   move equivalent to the chart $\Gamma_0$  together with a small copy of $\Gamma$ in $R_2$.  
A chart $\Gamma$ in a $2$-disk  is {\it nomadic} if it is nomadic with respect to every chart.  
}\end{definition}

\begin{figure}[h]
\begin{center}
\mbox{\epsfxsize=6cm \epsfbox{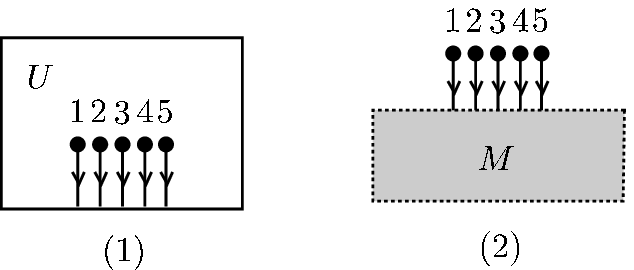}}
\end{center}
\vspace{-0.5cm}
\caption{}
\label{sfg09}
\end{figure}

\begin{lemma}\label{lem:nomadic}
Let $D$ be a $2$-disk and $B$ a compact, connected and oriented surface.  
\begin{itemize}
\item[{\rm (1)}] 
Let $\Gamma$ be a chart in $D$.  
If there is a $2$-disk $U$ in $D$ such that $\Gamma \cap U$ is as in Figure~$\ref{sfg09} (1)$, then
$\Gamma$ is nomadic.  
\item[{\rm (2)}] 
Let $\Gamma_0$ be a chart in $B$. 
If there is a $2$-disk $U$ in $B$ such that $\Gamma_0 \cap U$ is as in Figure~$\ref{sfg09} (1)$, then 
any chart $\Gamma$ in a $2$-disk is nomadic with respect to $\Gamma_0$.  
\end{itemize}
\end{lemma}

\begin{figure}[h]
\begin{center}
\mbox{\epsfxsize=9.0cm \epsfbox{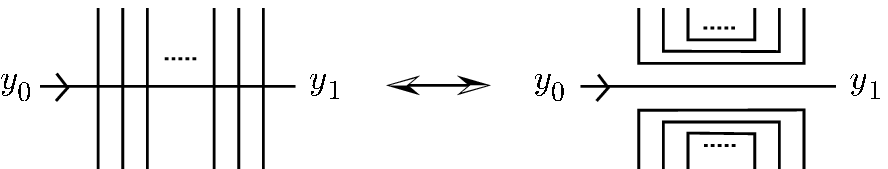}}
\end{center}
\vspace{-0.5cm}
\caption{}
\label{sfg10}
\end{figure}

\begin{figure}[h]
\begin{center}
\mbox{\epsfxsize=9.0cm \epsfbox{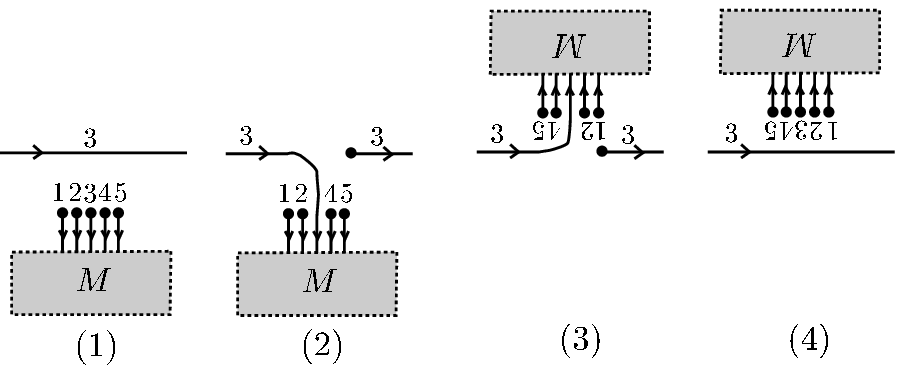}}
\end{center}
\vspace{-0.5cm}
\caption{}
\label{sfg11}
\end{figure}

\begin{figure}[h]
\begin{center}
\mbox{\epsfxsize=8.0cm \epsfbox{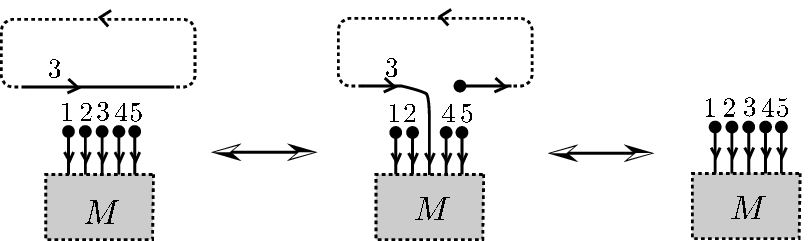}}
\end{center}
\vspace{-0.5cm}
\caption{}
\label{sfg12}
\end{figure}

{\it Proof.} 
(1) First we consider a special case where $\Gamma$ is as in Figure~$\ref{sfg09} (2)$.  
Let $\Gamma_0$ be any chart in $B$, and put a small copy of $\Gamma$ in a region of $B \setminus \Gamma_0$.  
As shown in Figure~\ref{sfg11}, it can pass through any edge of $\Gamma_0$ which is labeled in $\{1, \dots, 5\}$.  
For an edge labeled $\sigma$, apply a chart move as in Figure~$\ref{sfg04} (1)$, let  $\Gamma$ pass through the $12$ edges with labels $1$ and $2$, and recover the edge labeled $\sigma$ by the move in Figure~\ref{sfg04}.   Thus we see that $\Gamma$ is nomadic.  
Now we consider a general case.  Take a point $y_0$ in the region $U$ and a point $y_1$ in the boundary $\partial D$.  Consider a simple path $\eta : [0,1] \to D$ connecting $y_0$ and $y_1$ such that $\eta$ intersects $\Gamma$ transversely.  Let $w$ be the intersection word of $\eta$ with respect to $\Gamma$ (see \cite{Kam02, Kam07}).  
Let $\Gamma'$ be a chart obtained from $\Gamma$ by adding some hoops surrounding $\Gamma$ such that 
the intersection word $w'$ of $\eta$ with respect to $\Gamma'$ is $w \cdot w^{-1}$.   Applying a chart move in a neighborhood of $\eta$ as in Figure~\ref{sfg10}, 
we have a chart $\Gamma''$ such that it coincides with $\Gamma'$ outside of the neighborhood of $\eta$ and the path $\eta$ misses $\Gamma''$. 
So $\Gamma''$ is as in  Figure~$\ref{sfg09} (2)$.  Note that $\Gamma''$ is chart move equivalent to $\Gamma$, since 
one can add or remove any hoop surrounding it by chart moves as in Figure~$\ref{sfg12}$.  Since $\Gamma''$ is nomadic as shown in the previous case, we see that $\Gamma$ is nomadic.  

\begin{figure}[h]
\begin{center}
\mbox{\epsfxsize=8.5cm \epsfbox{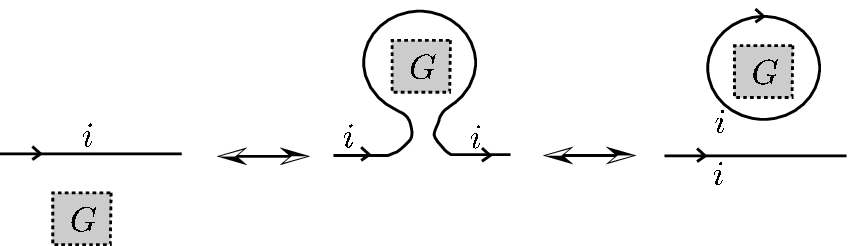}}
\end{center}
\vspace{-0.5cm}
\caption{}
\label{sfg13}
\end{figure}

Now we prove (2). Let $U$ be a region such that $\Gamma_0 \cap U$ is as in Figure~$\ref{sfg09} (1)$.  
It is sufficient to show that any chart $\Gamma$ put in a region of $B \setminus \Gamma_0$ can be moved into $U$.  
As shown in Figure~$\ref{sfg13}$,  $\Gamma$ can pass through any edge of $\Gamma_0$ by getting a surrounding hoop.   When $\Gamma$ arrives in $U$, it is surrounded some hoops, which can be removed by use of the edges of $\Gamma_0$ in $U$ as in Figure~$\ref{sfg12}$.  \qed 

\begin{figure}[h]
\begin{center}
\mbox{\epsfxsize=9.0cm \epsfbox{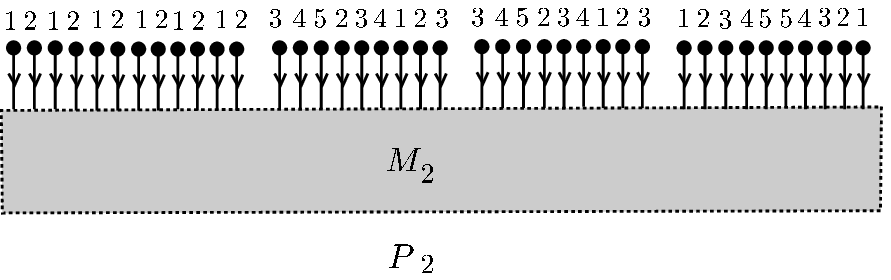}}
\end{center}
\vspace{-0.5cm}
\caption{Chart $P_2$, which is equivalent to $2 N_0$}
\label{sfg14}
\end{figure}

\begin{lemma}\label{lem:2n0}
Let $P_2$ be a chart depicted in Figure~$\ref{sfg14}$, where $M_2$ is the chart depicted in $Figure~\ref{sfg06}$.  It is chart move equivalent to $2 N_0$.  
\end{lemma}

\begin{figure}[h]
\begin{center}
\mbox{\epsfxsize=7.0cm \epsfbox{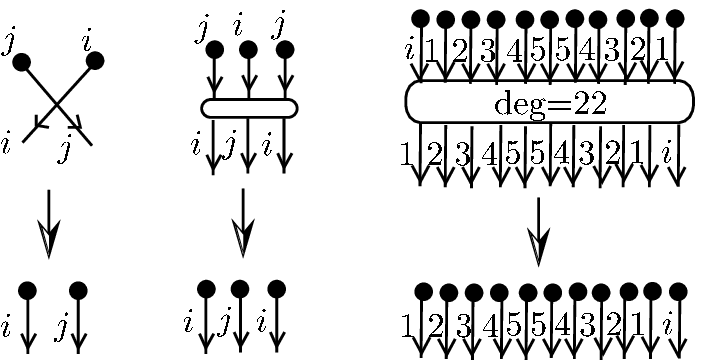}}
\end{center}
\vspace{-0.5cm}
\caption{Chart moves}
\label{sfg15}
\end{figure}

{\it Proof.} Applying chart moves depicted in Figure~\ref{sfg15} to the chart $P_2$, we have $2 N_0$.  \qed 

Now we prove Theorem~\ref{thm:chart1}.  

\begin{figure}[h]
\begin{center}
\mbox{\epsfxsize=8.5cm \epsfbox{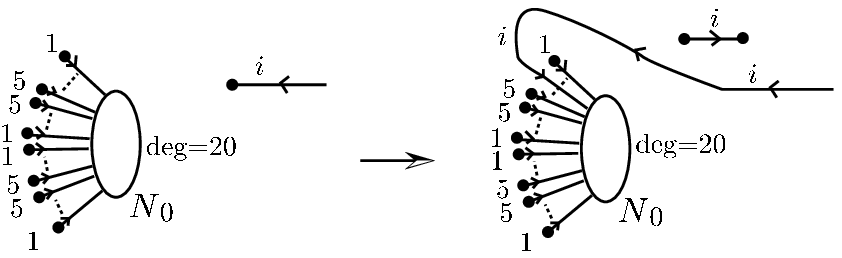}}
\end{center}
\vspace{-0.5cm}
\caption{}
\label{sfg16}
\end{figure}

\begin{figure}[h]
\begin{center}
\mbox{\epsfxsize=6.5cm \epsfbox{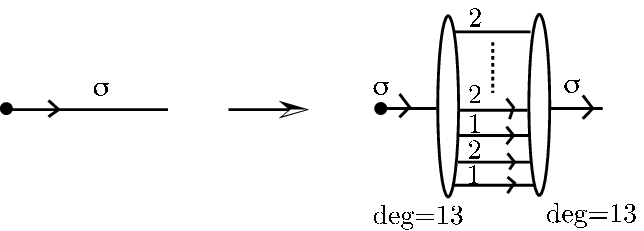}}
\end{center}
\vspace{-0.5cm}
\caption{}
\label{sfg17}
\end{figure}

{\it Proof of Theorem~\ref{thm:chart1}.}  
First we consider a case where $\Gamma$ is a chart with $n_{\rm II}^-(\Gamma)=0$.  
It suffices to show that 
$\Gamma \oplus (n_{\rm I}^-(\Gamma)  +   2 n_{\rm II}^+(\Gamma) +1)   \, N_0$ is chart move equivalent to 
$$\Gamma'  \oplus n_{\rm II}^+(\Gamma) \, N_2  \oplus n_{\rm I}^-(\Gamma) \, F_1  
$$
for some chart $\Gamma'$ with $n_{\rm I}^-(\Gamma')=n_{\rm II}^+(\Gamma') = n_{\rm II}^-(\Gamma')=0$ 
such that $\Gamma'$ has $N_0$ as a factor.  
By Lemma~\ref{lem:nomadic}, $N_0$ is nomadic.  Thus we can move $N_0$ freely up to chart move equivalence.  
For each black vertex of type ${\rm I}^-$, move a chart $N_0$ near the vertex and apply a chart move as in Figure~\ref{sfg16} to make a free edge.  Move the free edge toward the boundary of $B$ by the chart move as in Figure~\ref{sfg13}.  
Since there is at least one $N_0$ near $\partial B$, the hoops surrounding the free edge can be removed (Figure~\ref{sfg12}), and we may also assume that the label of the free edge is $1$ (Lemma~18.24 of \cite{Kam02}).  
Thus we can change $\Gamma \oplus (n_{\rm I}^-(\Gamma)  +   2 n_{\rm II}^+(\Gamma) +1)   \, N_0$  so that 
all black vertices of type ${\rm I}^-$ are endpoints of $F_1$'s near $\partial B$.  
We still have $2 n_{\rm II}^+(\Gamma) +1$ $N_0$'s near $\partial B$.  For each black vertex of type ${\rm II}^+$, move a pair of $N_0$ near the vertex.   Change the pair of $N_0$'s to a chart $P_2$ in Figure~$\ref{sfg14}$ (Lemma~\ref{lem:2n0}).  
The edge adjacent to the vertex of type ${\rm II}^+$ is oriented outward and is labeled $\sigma$.  Apply a chart move as in Figure~\ref{sfg17}, and then apply a chart move between the $12$ edges there and the $12$ edges of $P_2$ to get one $N_2$.  Move the chart $N_2$ toward $\partial B$.  (Note that $N_2$ is nomadic by Lemma~\ref{lem:nomadic}.)  Now all black vertices of type ${\rm II}^+$ belong to $N_2$'s near $\partial B$.  
We still have one $N_0$ near $\partial B$.  Thus the chart is 
$\Gamma'  \oplus n_{\rm II}^+(\Gamma) \, N_2  \oplus n_{\rm I}^-(\Gamma) \, F_1  
$ for a chart $\Gamma'$ with $n_{\rm I}^-(\Gamma')=n_{\rm II}^+(\Gamma') = n_{\rm II}^-(\Gamma')=0$ 
such that $\Gamma'$ has $N_0$ as a factor.  

We consider a case where $\Gamma$ is a chart with $n_{\rm II}^+(\Gamma) \geq n_{\rm II}^-(\Gamma)>0$.  
Let $v$ be a black vertex of type ${\rm II}^-$.  Choose a black vertex $v'$ of type ${\rm II}^+$ and consider a simple path $\eta$ from $v$ to $v'$.  If $\eta$ intersects an edge labeled $\sigma$, then apply a chart move depicted in Figure~$\ref{sfg04} (1)$ and we assume that $\eta$ intersects only edges with labels in $\{1, \dots, 5\}$.  For each intersection of $\eta$ and the chart, we assert one $N_0$ and apply a chart move as in Figure~$\ref{sfg11} (2)$ so that $\eta$ does not intersect the chart.  Now move $v$ along $\eta$ toward $v'$ and by a chart move we can make a free edge with label $\sigma$, that is $F_2$.  Move this $F_2$ toward $\partial B$ by 
moves as in Figure~\ref{sfg13}.   The hoops surrounding the free edge can be removed by adding one $N_0$ near $\partial B$ as before.  By this procedure, we can move all black vertices of type ${\rm II}^-$ near $\partial B$ as endpoints of $F_2$'s.  
The number of $F_2$'s is $n_{\rm II}^-(\Gamma)$.  
There are $n_{\rm II}^+(\Gamma) - n_{\rm II}^-(\Gamma)$ 
black vertices of type ${\rm II}^+$ in the chart, besides the endpoints of $F_2$'s.  For each black vertex of type ${\rm II}^+$, that is not an endpoint of $F_2$, add a pair of $N_0$ to make $P_2$.  As in the previous case, we can move the black vertex of type ${\rm II}^+$ as an endpoint of $N_2$ near $\partial B$.  The number of $N_2$'s is $n_{\rm II}^+(\Gamma) - n_{\rm II}^-(\Gamma)$.  As in the previous case, we move black vertices of type ${\rm I}^-$ as endpoints of $F_1$'as near $\partial B$.  The number of $F_1$'s is $n_{\rm I}^-(\Gamma)$.   Thus we have a chart written as 
$$\Gamma'  \oplus (n_{\rm II}^+(\Gamma)-  n_{\rm II}^-(\Gamma)) N_2  \oplus n_{\rm I}^-(\Gamma) \, F_1  
\oplus  n_{\rm II}^-(\Gamma) F_2  $$
for some chart $\Gamma'$ with $n_{\rm I}^-(\Gamma')=n_{\rm II}^+(\Gamma') = n_{\rm II}^-(\Gamma')=0$ 
such that $\Gamma'$ has $N_0$ as a factor.  \qed


\section{Proof of Theorem~\ref{thm:main}}\label{sect:mainproof}

{\it Proof of Theorem~\ref{thm:main}.}  
By Theorem A of Siebert and Tian \cite{ST03} (cf. \cite{Auroux2005}) a chiral and irreducible genus-two Lefschetz fibration $f'$ 
with transitive monodromy representation is holomorphic and hence it is a fiber sum of some copies of $f_0$ and $f_1$.  
Therefore Corollary~\ref{cor:chart1} implies the assertions (2) and (4) of Theorem~\ref{thm:main}, and the former part of (3).    
It is well-known that $3 f_0 \cong 2 f_1$.   Thus we can take $b$ to be $0$ or $1$.  
We shall compare the number of singular fibers of each type of $f \,\# \, m   \,  f_0 $ with that of 
 ${\#}  \, (a+m)  \,  f_0 
 \,  \# \,  b  \,  f_1 
 \,  \#  \, c   \,  f_2
 \,  \#  \, d   \, f'_1
 \,  \#  \,  e  \,  f'_2 $. 
We have already used the information on the numbers of singular fibers of type ${\rm I}^-$, ${\rm II}^+$ and type ${\rm II}^-$ to determine $c$, $d$ and $e$; $c = n_{\rm II}^+(f) - n_{\rm II}^-(f)$,  
$d = n_{\rm I}^-(f)$ and $e = n_{\rm II}^-(f)$.   
The number of singular fibers of type ${\rm I}^+$ of 
$f \,\# \, m   \,  f_0 $ is $n_{\rm I}^+(f) + 20 m$, and that of 
${\#}  \, (a+m)  \,  f_0 
 \,  \# \,  b  \,  f_1 
 \,  \#  \, c   \,  f_2
 \,  \#  \, d   \, f'_1
 \,  \#  \,  e  \,  f'_2 $
is $20 (a+m) + 30 b + 28c + d$.  
From this equality, we have 
$$20 a + 30 b = n_{\rm I}^+(f) - n_{\rm I}^-(f) 
-28 (n_{\rm II}^+(f)  
- n_{\rm II}^-(f)).$$
Thus the right hand side, which is ${\cal E}(f)$, is a multiple of $10$.  
And we see that the parity of $b$ equals to the parity $\epsilon(f)$.   
Therefore when we  assume $b=\epsilon(f)$, we have $a = ({\cal E}(f) -30 \epsilon(f))/20$.  \qed

\section{Concluding remark}\label{sect:concluding}

In the proof of Theorem~\ref{thm:main}, we assumed the deep result due to Siebert and Tian \cite{ST03} stating that 
any chiral and irreducible genus-two Lefschetz fibration 
$f'$ over $S^2$ with transitive monodromy representation 
is holomorphic and it is a fiber sum of some copies of $f_0$ and $f_1$.  If  one does not need a lower bound $m_0$ given in (4) of Theorem~\ref{thm:main}, we can prove Theorem~\ref{thm:main} without assuming  Siebert and Tian's result.  

\begin{proposition}\label{prop:chart2}
Let $\Gamma$ be a chart description of a chiral and irreducible genus-two Lefschetz fibration over $B=D^2$ (or $S^2$).  
There exists a positive integer $m$ such that 
$\Gamma \oplus m N_0$ is chart move equivalent to $ (a + m) N_0 \oplus b N_1$ for some 
integers $a$ and $b$.  
\end{proposition}

\begin{figure}[h]
\begin{center}
\mbox{\epsfxsize=12.5cm \epsfbox{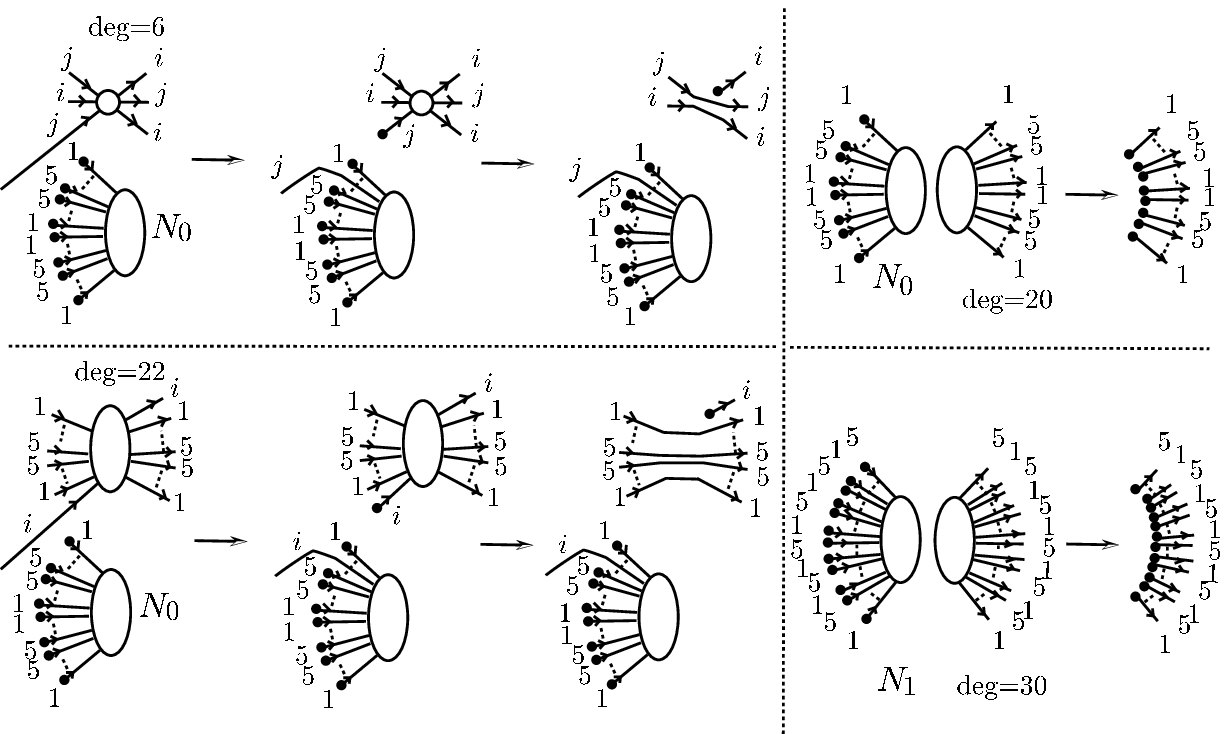}}
\end{center}
\vspace{-0.5cm}
\caption{}
\label{sfg18}
\end{figure}

{\it Proof.}  Since $f$ is chiral and irreducible,  we may assume that $\Gamma$ is chiral and irreducible by 
Proposition~\ref{prop:chartdescriptionchiralirreducible}.  
Adding some $N_0$'s to the chart and applying chart moves shown in Figure~\ref{sfg18},  
we can remove all degree-$6$ vertices, degree-$22$ vertices,  degree-$20$ vertices whose adjacent edges are oriented outward, and 
degree-$30$ vertices whose adjacent edges are oriented outward.    (Since $3 N_0$   is chart move equivalent to  $2 N_1$, we may add $N_1$'s too.)  
Remove all hoops using an $N_0$ (Figure~$\ref{sfg12}$).  
Now every edge is adjacent to a black vertex, a degree-$4$ vertex, a degree-$20$ vertex whose adjacent edges are oriented inward or a degree-$30$ vertex whose adjacent edges are oriented inward.   Note that for a degree-$4$ vertex, the two incoming adjacent edges have black vertices at the other end.  Thus by a chart move depicted in Figure~\ref{sfg15}, we can remove the degree-$4$ vertex.  Remove all degree-$4$ vertices this way.  Now the chart is a union of some $N_0$'s and $N_1$'s.  \qed

Now we have a corollary to Proposition~\ref{prop:chart2}.  

\begin{corollary}
Let $f$ be a chiral and irreducible genus-two Lefschetz fibration over $S^2$.  
There exists a positive number $m$ such that 
$f \#  m f_0 \cong (a + m) f_0 \# b f_1$ for some 
integers $a$ and $b$.  
\end{corollary} 

Using this corollary, we have a proof of Theorem~\ref{thm:main}, except the assertion (4), without using Siebert and Tian's result.  



\begin{thebibliography}{99}

\bibitem{Auroux2003} 
D. Auroux, {\it Fiber sums of genus $2$ Lefschetz fibrations}, 
Turkish J. Math. {\bf 27} (2003) 1--10. 

\bibitem{Auroux2005}
D. Auroux, {\it A stable classification of Lefschetz fibrations}, 
Geom. Topol. {\bf 9} (2005), 203--217. 

\bibitem{Birman}
J. Birman, {\it Braids, links and mapping class groups}, Princeton Univ. Press, 1974. 

\bibitem{GS} 
R. E. Gompf and A. I. Stipsicz, 
{\it $4$-manifolds and Kirby calculus}, 
Graduate Studies in Math. {\bf 20}, 
Amer. Math. Soc., Providence, RI, 1999. 


\bibitem{Kam92}
S. Kamada, 
{\it Surfaces in $R^4$ of braid index three are ribbon}, 
J. Knot Theory Ramifications {\bf 1} (1992) 137--160.

\bibitem{Kam96}
S. Kamada, 
{\it An observation of surface braids via chart description},  
J. Knot Theory Ramifications {\bf 4} (1996), 517--529. 

\bibitem{Kam02}
S. Kamada, 
{\it Braid and knot theory in dimension four}, 
Math. Surveys Monogr. {\bf 95}, 
Amer. Math. Soc., Providence, RI, 2002.

\bibitem{Kam07} 
S. Kamada, 
{\it Graphic descriptions of monodromy representations}, 
Topology Appl. {\bf 154}  (2007),  1430--1446.

\bibitem{KMMW} 
S. Kamada, Y. Matsumoto, T. Matumoto and K. Waki, 
{\it Chart description and a new proof of the classification theorem of genus
one  Lefschetz fibrations}, 
J. Math. Soc. Japan,  {\bf 57} (2005), 537--555. 

\bibitem{Kas80}
A. Kas, 
{\it On the handlebody decomposition associated to a Lefschetz
fibration}, 
Pacific J. Math.  {\bf 89} (1980), 89--104. 


\bibitem{MH}
R. Mandelbaum and J. R. Harper, 
{\it Global monodromy of elliptic Lefschetz fibrations}, 
in ``Current trends in algebraic topology'', 
pp. 35--41, CMS Conf. Proc. {\bf 2}, Amer. Math. Soc., Providence, RI, 1982. 




\bibitem{Matsu86}
Y. Matsumoto, 
{\it Diffeomorphism types of elliptic surfaces}, 
Topology  {\bf 25} (1986), 549--563. 

\bibitem{Matsu96}
Y. Matsumoto, 
{\it  Lefschetz fibrations of genus two
 - A topological approach},
in ``Topology and Teichm\"uller spaces'' 
S. Kojima et al, eds., Proc. the 37-th Taniguchi Sympo.,  pp. 123--148,  
World Scientific Publishing, River Edge, NJ, 1996.

\bibitem{Moi77}
B. G. Moishezon, 
{\it Complex surfaces and connected sums of complex projective planes}, 
Lecture Notes in Math. {\bf 603}, Springer Verlag, 1977. 

\bibitem{Moi81}
B. G. Moishezon, 
{\it Stable branch curves and braid monodromies}, 
in ``Algebraic Geometry'', pp. 107--192, 
Lecture Notes in Math. {\bf 862}, Springer Verlag, 1981. 

\bibitem{ST99} 
B. Siebert and G. Tian, 
{\it On hyperelliptic $C^\infty$-Lefschetz fibrations of four-manifolds}, 
Commun. Contemp. Math. {\bf 1} (1999), no. 2, 255--280. 


\bibitem{ST03}
B. Siebert and G. Tian, 
{\it On the holomorphicity of genus two Lefschetz fibrations}, 
Ann. of Math. (2) {\bf 161} (2005), no. 2, 959--1020.

\end{thebibliography}
\end{document}